\documentclass[11pt,oneside]{article}
\usepackage{mathrsfs}
\usepackage{amssymb}
\usepackage{amsfonts}
\usepackage{bbding}
\usepackage{pifont}
\usepackage{fancyhdr}
\usepackage{titlesec}
\usepackage{cite}
\usepackage{ifthen}
\usepackage{amssymb,amsthm,amscd}

\titleformat{\section}{\centering\large\bfseries}{\S\arabic{section}}{1em}{}
\newboolean{first}
\setboolean{first}{true}

\begin{document}

\setlength\abovedisplayskip{2pt}
\setlength\abovedisplayshortskip{0pt}
\setlength\belowdisplayskip{2pt}
\setlength\belowdisplayshortskip{0pt}

\title{\bf \huge Symmetries of geometric flows}
\author{Xu Chao\footnote{e-mail address: xuchaomykj@163.com}\\
        Department of Mathematics\\
        Zhejiang University, Hangzhou, China} \maketitle

\begin{abstract}
By applying the theory of group-invariant solutions we investigate
the symmetries of Ricci flow and hyperbolic geometric flow both on
Riemann surfaces. The warped products on $\mathcal {S}^{n+1}$ of
both flows are also studied.\\\hskip 2mm {\small }
\end{abstract}

\thispagestyle{fancyplain} \fancyhead{}

\section{Introduction}

The Ricci flow is the geometric evolution equation in which one
starts with a smooth Riemannian manifold $(\mathcal {M}^n,g_0)$ and
evolves its metric by the equation
\\$$\frac{\partial}{\partial t}g=-2Rc,\eqno(1.1)$$\\
where $Rc$ denotes the Ricci tensor of the metric $g$. The Ricci
flow has been exhaustively studied and successfully applied to solve
the famous Poincar$\acute{e}$'s Conjecture \cite{CZ}. Recently, Kong
and Liu \cite{KL} introduced the hyperbolic geometric flow which is
the hyperbolic version of Ricci flow
\\$$\frac{\partial^2}{\partial t^2}g=-2Rc,\eqno(1.2)$$\\
which shows different behavior with the original Ricci flow.

On Riemann surfaces $(\mathcal{M}^2,g)$, equations (1.1) and (1.2)
can be simplified to scalar equations
\\$$u_t=\triangle\ln u,\eqno(1.3)$$
$$u_{tt}=\triangle\ln u,\eqno(1.4)$$\\
where function $u(x,y,t)$ is the conformal factor of $g$:
\\$$g_{ij}=u(x,y,t)\delta_{ij}$$\\
Later we will use the theory of group-invariant solutions to
investigate (1.3) and (1.4). As we will see, the sets of symmetries
of the two equations are quite large and we expect to find large
classes of exact solutions to both flows on Riemann surfaces and the
symmetries dependent on the solution of two-dimensional Laplace
equation.

The technique here we use to investigate the symmetries and exact
solutions of the equations is the theory of group-invariant
solutions for differential equations which applies Lie group, Lie
algebra and adjoint representation to differential equations. For
most cases, there is a one-to-one correspondence between different
symmetries of an equation and the conjugate classes of subgroups of
its one-parameter transformation group. So finally through the
classification of subalgebras of the Lie algebra of the
transformation group, we are able to classify all the symmetries of
the equations. We will introduce this technique briefly later in
this paper. For more details, see \cite{O}.

We will also study warped products on $\mathcal {S}^{n+1}$ or
$SO(n+1)$-invariant metrics on $\mathcal {S}^{n+1}$ of both flows on
the set $(-1,1)\times\mathcal {S}^n$:
\\$$g=\varphi^2(x,t)dx^2+\psi^2(x,t)g_{can},$$\\
where $g_{can}$ denotes the canonical metric on $\mathcal {S}^n$.
This metric under Ricci flow was studied in [1]. Analyzing its
asymptotic behavior leads to significant information about the
nechpinch which is important in the surgery of Ricci flow. By a
change of coordinate
\\$$s(x)=\int^x_0\varphi(x)dx,$$\\
the evolutions of $\varphi(s,t)$ and $\psi(s,t)$ under Ricci flow
and hyperbolic geometric flow are the followings respectively:
\\$$\left\{\begin{array}{ll}
\varphi_t=n\frac{\psi_{ss}}{\psi}\varphi\\
\psi_t=\psi_{ss}-(n-1)\frac{1-\psi^2_s}{\psi}
\end{array}\right.\eqno(1.5)$$\\
under Ricci flow, and
\\$$\left\{\begin{array}{ll}
\varphi_{tt}=n\frac{\psi_{ss}}{\psi}\varphi-\frac{\varphi^2_t}{\varphi}\\
\psi_{tt}=\psi_{ss}-(n-1)\frac{1-\psi^2_s}{\psi}-\frac{\psi^2_t}{\psi}
\end{array}\right.\eqno(1.6)$$\\
under hyperbolic geometric flow. In contrast to (1.3) and (1.4),
equations (1.5) and (1.6) have few symmetries especially in higher
dimensions.

This paper is organized as follows: We would begin with the theory
of group-invariant solutions for differential equations in Section
2. In Section 3, we will study the symmetries and exact solutions of
Ricci flow on surfaces. In Section 4, we investigate hyperbolic
geometric flow on surfaces. In section 5, the warped product of
$\mathcal {S}^{n+1}$ on both flows are studied. In Section 6, we
give some further discussions. Finally, in section 7, we derive the
evolutions
of warped products on both flow.\\

{\bf Acknowledgement. } {\it The author thanks the Center of
Mathematical Sciences at Zhejiang University where he wrote this
paper during the summer of 2009. }\\


\section{Theory of group-invariant solutions for differential
equations}

In this section, we briefly introduce the theory of group-invariant
solutions for differential equations. The following main definitions
and theorems are cited from \cite{O}.

First we introduce the jet space. Given
\\$$u_t=\triangle\ln u,$$\\
let $w=\ln u$, so
\\$$e^ww_t-w_{xx}-w_{yy}=0.\eqno(2.1)$$\\
We regard $w$ and its derivatives as variables in (2.1), so (2.1)
can be regarded as defined on
\\$$X\times U^{(2)}=\{(x,y,t;w;w_x,w_y,w_t;w_{xx},w_{xy},w_{xt},w_{yy},w_{yt},w_{tt})\},$$\\
where $X=\{(x,y,t)\}$ is the space of independent variables
$(x,y,t)$.

In general, we denote an $n$-th order differential equation of $w$
with independent variables $x=(x^1,...,x^p)$ by
\\$$\triangle (x,w^{(n)})=0.$$\\
Thus $\triangle$ can be regarded as a smooth map form the {\it jet
space} $X\times U^{(n)}$ to $\mathbb{R}$
\\$$\triangle:X\times U^{(n)}\rightarrow \mathbb{R},$$\\
and the differential equation tells  where the given map $\triangle$
vanishes on $X\times U^{(n)}$, thus determines a subvariety
\\$$\mathscr{S}_{\triangle}=\{(x,w^{(n)}):\triangle (x,w^{(n)})=0\}\subset X\times U^{(n)}$$\\
of the total jet space.\\

{\bf Definition 2.1.} {\it Let $\mathscr{S}$ be a system of
differential equations. A symmetry group of the system $\mathscr{S}$
is a local group of transformations $G$ acting on an open subset $M$
of the space of independent and dependent variables for the system
with the property that whenever $u=f(x)$ is a solution of
$\mathscr{S}$, and whenever $g\cdot f$ is defined for $g\in G$, then
$u=g\cdot f(x)$ is also a solution of the system.}\\

{\bf Theorem 2.2.} {\it Let $M$ be an open subset of $X\times
U^{(n)}$ and suppose $\triangle (x,w^{(n)})=0$ is an $n$-th order
equation defined over $M$, with corresponding subvariety
$\mathscr{S}_{\triangle}\subset M$. Suppose $G$ is a local group of
transformations acting on $M$ which leaves $\mathscr{S}_{\triangle}$
invariant, meaning that whenever
$(x,w^{(n)})\in\mathscr{S}_{\triangle}$, we have $g\cdot
(x,w^{(n)})\in\mathscr{S}_{\triangle}$ for all $g\in G$ such that
this is defined. Then $G$ is a symmetry group of the equation in the
sense of Definition 2.1.}\\

Next we introduce the prolongation of vector fields corresponding to
one-parameter transformation group acting on $M\subset X\times
U=\{(x,w)\}$. We only state its formula here for our use. The
interested reader can see \cite{O}.\\

{\bf Theorem 2.3.} {\it Let
\\$${\bf v}=\sum_{i=1}^p\xi^i(x,w)\frac{\partial}{\partial x^i}+\phi(x,w)\frac{\partial}{\partial w}$$\\
be a vector field defined on an open subset $M\subset X\times U$.
The $n$-th prolongation of ${\bf v}$ is the vector field
\\$$pr^{(n)}{\bf v}={\bf v}+\sum_J\phi^J(x,w^{(n)})\frac{\partial}{\partial w_J}$$\\
defined on the corresponding jet space $M^{(n)}\subset X\times
U^{(n)}$, the summation being over all (unordered) multi-indices
$J=(j_1,...,j_k)$, with $1\leq j_k\leq p$, $1\leq k\leq n$. The
coefficient functions $\phi^J$ of $pr^{(n)}{\bf v}$ are given by the
following formula:
\\$$\phi^J(x,w^{(n)})=D_J(\phi-\sum_{i=1}^p\xi^iw_i)+\sum_{i=1}^p\xi^iw_{J,i},$$\\
where $D$ is the total derivative operator, and $w_i=\partial
w/\partial x^i$, $w_{J,i}=\partial
w_J/\partial x^i$.}\\

We state two important definitions which play important role in the
theory.\\

{\bf Definition 2.4.} {\it Let
\\$$\triangle (x,w^{(n)})=0,$$\\
be a differential equation. The equation is said to be of maximal
rank if the Jacobian matrix
\\$$\mathcal {J}_{\triangle}(x,w^{(n)})=(\frac{\partial\triangle}{\partial x^i},\frac{\partial\triangle}{\partial w_J})$$\\
of $\triangle$ with respect to all the variables $(x,w^{(n)})$ is of
rank $1$ whenever $\triangle (x,w^{(n)})=0$.}\\

For example, consider (2.1), the corresponding Jacobian matrix is
\\$$\mathcal {J}_{\triangle}(x,y,t;w;w_x,w_y,w_t;w_{xx},w_{xy},w_{xt},w_{yy},w_{yt},w_{tt})
=(0,0,0;e^ww_t;0,0,e^w;-1,0,0,-1,0,0),$$\\
which is of rank $1$ whenever $\triangle (x,y,t,w^{(2)})=0$. So
(2.1) is of maximal rank.\\

{\bf Definition 2.5.} {\it An $n$-th order differential equation
$\triangle (x,w^{(n)})=0$ is locally solvable at the point
\\$$(x_0,w_0^{(n)})\in\mathscr{S}_{\triangle}=\{(x,w^{(n)}):\triangle (x,w^{(n)})=0\}$$\\
if there exists a smooth solution $u=f(x)$ of the equation, defined
for $x$ in a neighborhood of $x_0$, which has the prescribed initial
condition $w_0^{(n)}=pr^{(n)}f(x_0)$, where $pr^{(n)}f(x_0)$ means
$f$ and all its derivatives up to order $n$ at point $x_0$. The
equation is locally solvable if it is locally solvable at every
point of $\mathscr{S}_{\triangle}$. A differential equation is
nondegenerate if at every point
$(x_0,w_0^{(n)})\in\mathscr{S}_{\triangle}$ it is both locally
solvable and of maximal rank.}\\

The main theorem we will use is the following:\\

{\bf Theorem 2.6.} {\it Let $\triangle (x,w^{(n)})=0$ be a
nondegenerate differential equation. A connected local group of
transformations $G$ acting on an open subset $M\subset X\times U$ is
a symmetry group of the equation if and only if
\\$$pr^{(n)}{\bf v}[\triangle (x,w^{(n)})]=0,\qquad whenever\qquad\triangle(x,w^{(n)})=0,$$\\
for every infinitesimal generator {\bf v} of $G$.}\\

We calculate some prolongation formulas here that we will use later.
On $M\subset X\times U$, given a vector
\\$${\bf v_1}=\xi(x,y,t,w)\frac{\partial}{\partial x}+\eta(x,y,t,w)\frac{\partial}{\partial y}
+\tau(x,y,t,w)\frac{\partial}{\partial t}+\phi(x,y,t,w)\frac{\partial}{\partial w},$$\\
its second order prolongation is
\\$$pr^{(2)}{\bf v_1}={\bf v_1}+\phi^x\frac{\partial}{\partial w_x}+\phi^y\frac{\partial}{\partial w_y}
+\phi^t\frac{\partial}{\partial
w_t}+\phi^{xx}\frac{\partial}{\partial
w_{xx}}+\phi^{xy}\frac{\partial}{\partial w_{xy}}
+\phi^{xt}\frac{\partial}{\partial
w_{xt}}+\phi^{yy}\frac{\partial}{\partial w_{yy}}
+\phi^{yt}\frac{\partial}{\partial w_{yt}}+\phi^{tt}\frac{\partial}{\partial w_{tt}}.$$\\
We will use the followings:
\\$$\phi^t=\phi_t-\xi_tw_x-\eta_tw_y+(\phi_w-\tau_t)w_t-\xi_ww_xw_t-\eta_ww_yw_t-\tau_ww_t^2,$$
\\$$\phi^x=\phi_x+(\phi_w-\xi_x)w_x-\eta_xw_y-\tau_xw_t-\xi_ww_x^2-\eta_ww_xw_y-\tau_ww_xw_t,$$
\\$$\phi^y=\phi_y-\xi_yw_x+(\phi_w-\eta_y)w_y-\tau_yw_t-\xi_ww_xw_y-\eta_ww_y^2-\tau_ww_yw_t,$$
\\$$\begin{array}{ll}
\phi^{tt}&=\phi_{tt}+(2\phi_{tw}-\tau_{tt})w_t-\eta_{tt}w_y-\xi_{tt}w_x+(\phi_{ww}-2\tau_{tw})w_t^2-2\eta_{tw}w_yw_t\cr
&-2\xi_{tw}w_xw_t-\tau_{ww}w_t^3-\eta_{ww}w_t^2w_y-\xi_{ww}w_t^2w_x+(\phi_w-2\tau_t)w_{tt}\cr
&-2\xi_tw_{xt}-w\eta_tw_{yt}-3\tau_ww_tw_{tt}-\eta_ww_yw_{tt}-\xi_ww_xw_{tt}\cr
&-2\eta_ww_tw_{yt}-2\xi_ww_tw_{xt},
\end{array}$$
\\$$\begin{array}{ll}
\phi^{xx}&=\phi_{xx}+(2\phi_{xw}-\xi_{xx})w_x-\eta_{xx}w_y-\tau_{xx}w_t+(\phi_{ww}-2\xi_{xw})w_x^2-2\eta_{xw}w_xw_y\cr
&-2\tau_{xw}w_xw_t-\xi_{ww}w_x^3-\eta_{ww}w_x^2w_y-\tau_{ww}w_x^2w_t+(\phi_w-2\xi_x)w_{xx}\cr
&-2\tau_xw_{xt}-2\eta_xw_{xy}-3\xi_ww_xw_{xx}-\eta_ww_yw_{xx}-\tau_ww_tw_{xx}\cr
&-2\eta_ww_xw_{xy}-2\tau_ww_xw_{xt},
\end{array}$$
\\$$\begin{array}{ll}
\phi^{yy}&=\phi_{yy}+(2\phi_{yw}-\eta_{yy})w_y-\xi_{yy}w_x-\tau_{yy}w_t+(\phi_{ww}-2\eta_{yw})w_y^2-2\xi_{yw}w_xw_y\cr
&-2\tau_{yw}w_yw_t-\eta_{ww}w_y^3-\xi_{ww}w_y^2w_x-\tau_{ww}w_y^2w_t+(\phi_w-2\eta_y)w_{yy}\cr
&-2\tau_yw_{yt}-2\xi_yw_{xy}-3\eta_ww_yw_yy-\xi_ww_xw_yy-\tau_ww_tw_{yy}\cr
&-2\xi_ww_yw_{xy}-2\tau_ww_yw_{yt}.
\end{array}$$\\

Next, given
\\$${\bf v_2}=\xi(s,t,\psi)\frac{\partial}{\partial s}+\tau(s,t,\psi)\frac{\partial}
{\partial t}+\phi(s,t,\psi)\frac{\partial}{\partial\psi},$$\\
its second order prolongation is
\\$$pr^{(2)}{\bf v_2}={\bf v_2}+\phi^s\frac{\partial}{\partial\psi_s}+\phi^t\frac{\partial}{\partial\psi_t}
+\phi^{ss}\frac{\partial}{\partial\psi_{ss}}+\phi^{st}\frac{\partial}{\partial\psi_{st}}
+\phi^{tt}\frac{\partial}{\partial\psi_{tt}}.$$\\

We will use the followings:
\\$$\phi^t=\phi_t-\xi_t\psi_s+(\phi_{\psi}-\tau_t)\psi_t-\xi_{\psi}\psi_s\psi_t-\tau_{\psi}\psi_t^2,$$
\\$$\phi^s=\phi_s+(\phi_{\psi}-\xi_s)\psi_s-\tau_s\psi_t-\xi_{\psi}\psi_s^2-\tau_{\psi}\psi_s\psi_t,$$
\\$$\begin{array}{ll}
\phi^{tt}&=\phi_{tt}+(2\phi_{t\psi}-\tau_{tt})\psi_t-\xi_{tt}\psi_s+(\phi_{\psi\psi}-2\tau_{t\psi})\psi_t^2
-2\xi_{t\psi}\psi_s\psi_t\cr
&-\tau_{\psi\psi}\psi_t^3-\xi_{\psi\psi}\psi_s\psi_t^2+(\phi_{\psi}-2\tau_t)\psi_{tt}-2\xi_t\psi_{st}
-3\tau_{\psi}\psi_t\psi_{tt}\cr
&-\xi_{\psi}\psi_s\psi_{tt}-2\xi_{\psi}\psi_t\psi_{st},
\end{array}$$
\\$$\begin{array}{ll}
\phi^{ss}&=\phi_{ss}+(2\phi_{s\psi}-\xi_{ss})\psi_s-\tau_{ss}\psi_t+(\phi_{\psi\psi}-2\xi_{s\psi})\psi_s^2
-2\tau_{s\psi}\psi_s\psi_t\cr
&-\xi_{\psi\psi}\psi_s^3-\tau_{\psi\psi}\psi_s^2\psi_t+(\phi_{\psi}-2\xi_s)\psi_{ss}-2\tau_s\psi_{st}
-3\xi_{\psi}\psi_s\psi_{ss}\cr
&-\tau_{\psi}\psi_t\psi_{ss}-2\tau_{\psi}\psi_s\psi_{st}.
\end{array}$$\\

When we solve out for example $\xi(x,y,t,w)$, $\eta(x,y,t,w)$,
$\tau(x,y,t,w)$, $\phi(x,y,t,w)$ for ${\bf v_1}$, we get some
vectors
\\$${\bf v_1},\cdot\cdot\cdot,{\bf v_k}.$$\\
These vectors generate the Lie algebra of the transformation group
$G$. To classify its subalgebras, we need to calculate the structure
constants
\\$$[{\bf v_i},{\bf v_j}]=C_{ij}^l{\bf v_l},$$\\
and the adjoint representations
\\$$\begin{array}{ll}
Ad(\exp(\varepsilon{\bf v_i})){\bf
v_j}&=\sum_{n=0}^{\infty}\frac{\varepsilon^n}{n!}(ad {\bf
v_i})^n({\bf v_j})\cr &={\bf v_j}-\varepsilon[{\bf v_i},{\bf
v_j}]+\frac{\varepsilon^2}{2}[{\bf v_i},[{\bf v_i},{\bf
v_j}]]-\cdot\cdot\cdot.
\end{array}$$\\

After the classification of subalgebras and subgroups of $G$, we get
an optimal system for the equation. By constructing invariants from
{\bf v} in these subalgebras, we can simplify the equation to ODE or
lower order PDE, thus we expect to find exact symmetric solutions to
the original equation. These will be investigated in detail for our
equations in the following sections.\\

\section{Ricci flow on Riemann surfaces}

On a surface, all of the information about curvature is contained in
the scalar curvature function $R$. The Ricci curvature is given by
\\$$R_{ij}=\frac{1}{2}Rg_{ij},$$\\
and the Ricci flow equation can be simplified to
\\$$\frac{\partial}{\partial t}g_{ij}=-Rg_{ij}.$$\\
The metric for a surface can always be written (at least locally) in
the following form
\\$$g_{ij}=u(x,y,t)\delta_{ij},$$\\
where $u(x,y,t)>0$. Therefore, we have
\\$$R=-\frac{\triangle\ln u}{u}.$$\\
Thus
\\$$\frac{\partial}{\partial t}u=\frac{\triangle\ln u}{u}\cdot u,$$\\
namely,
\\$$u_t-\triangle\ln u=0.\eqno(3.1)$$\\
Denote $w=\ln u$, thus
\\$$\triangle(x,y,t,w^{(2)})=w_{xx}-w_{yy}-e^ww_t=0\eqno(3.2)$$\\
we will use the techniques developed in section 2 to analyze (3.2).

First note that the Jacobian matrix of (3.2) is
\\$$\mathcal {J}_{\triangle}(x,y,t;w;w_x,w_y,w_t;w_{xx},w_{xy},w_{xt},w_{yy},w_{yt},w_{tt})
=(0,0,0;-e^ww_t;0,0,-e^w;1,0,0,1,0,0),$$\\
which is obviously of rank $1$ in $\mathscr{S}_{\triangle}$ and
(3.2) is obviously locally solvable. So we can apply Theorem 2.6.

Given a vector
\\$${\bf v}=\xi(x,y,t,w)\frac{\partial}{\partial x}+\eta(x,y,t,w)\frac{\partial}{\partial y}
+\tau(x,y,t,w)\frac{\partial}{\partial t}+\phi(x,y,t,w)\frac{\partial}{\partial w},$$\\
we have
\\$$pr^{(2)}{\bf v}(\triangle(x,y,t,w^{(2)}))=-\phi e^ww_t-\phi^te^w+\phi^{xx}+\phi^{yy}.$$\\
We apply the formulas for $\phi^t$, $\phi^{xx}$, $\phi^{yy}$ derived
in the above section, since $pr^{(2)}{\bf
v}(\triangle(x,y,t,w^{(2)}))=0$ whenever (3.2) holds, we use
$w_t=e^{-w}(w_{xx}+w_{yy})$ to cancel $w_t$ and set coefficients of
every monomial zero. For example the coefficient of $w_yw_{yt}$ is
$-2\tau_w$. So $\tau_w=0$, i.e. $\tau=\tau(x,y,t)$. See the
following table for all the coefficients.
\\$$\begin{tabular}{rr|rr}
\hline
monomial&coefficient&monomial&coefficient\\
$e^{-w}w_{xx}:$&$-\tau_{xx}-\tau_{yy}=0$&$e^{-w}w_{yy}:$&$-\tau_{xx}-\tau_{yy}=0$\\
$e^{-w}w_xw_{xx}:$&$-2\tau_{xw}=0$&$e^{-w}w_xw_{yy}:$&$-2\tau_{xw}=0$\\
$e^{-w}w_x^2w_{xx}:$&$-\tau_{ww}=0$&$e^{-w}w_x^2w_{yy}:$&$-\tau_{ww}=0$\\
$e^{-w}w_yw_{xx}:$&$-\tau_{yw}=0$&$e^{-w}w_yw_{yy}:$&$-2\tau_{yw}=0$\\
$e^{-w}w_y^2w_{xx}:$&$-\tau_{ww}=0$&$e^{-w}w_y^2w_{yy}:$&$-\tau_{ww}=0$\\
$e^w:$&$-\phi_t=0$&$e^ww_x:$&$\xi_t=0$\\
$e^ww_y:$&$\eta_t=0$&$1:$&$\phi_{xx}+\phi_{yy}=0$\\
$w_{xx}:$&$-\phi=\tau_t-2\xi_x=0$&$w_{yy}:$&$-\phi+\tau_t-2\eta_y=0$\\
$w_xw_{xx}:$&$-2\xi_w=0$&$w_yw_{yy}:$&$-2\eta_w=0$\\
$w_x:$&$2\phi_{xw}-\xi_{xx}-\xi_{yy}=0$&$w_y:$&$2\phi_{yw}-\eta_{xx}-\eta_{yy}=0$\\
$w_x^2:$&$\phi_{ww}-2\xi_{xw}=0$&$w_xw_y:$&$-2\eta_{xw}-2\xi_{yw}=0$\\
$w_x^3:$&$-\xi_{ww}=0$&$w_x^2w_y:$&$-\eta_{ww}=0$\\
$w_{xt}:$&$-\tau_x=0$&$w_{xy}:$&$-2\eta_x-2\xi_y=0$\\
$w_xw_{xy}:$&$-2\eta_w=0$&$w_xw_{xt}:$&$-2\tau_w=0$\\
$w_y^2:$&$\phi_{ww}-2\eta_{yw}=0$&$w_y^3:$&$-\eta_{ww}=0$\\
$w_y^2w_x:$&$-\xi_{ww}=0$&$w_{yt}:$&$-2\tau_y=0$\\
$w_yw_{xy}:$&$-2\xi_w=0$&$w_yw_{yt}:$&$-2\tau_w=0$\\
\hline
\end{tabular}$$\\

Thus we finally get
\\$$\left\{\begin{array}{ll}
\xi=\xi(x,y)\\
\eta=\eta(x,y)\\
\tau=c_1+c_2t\\
\phi=c_2-2\xi_x\\
\xi_x-\eta_y=0\\
\eta_x+\xi_y=0
\end{array}\right.\eqno(3.3)$$\\
From the last two relations in (3.3), we have
\\$$\left\{\begin{array}{ll}
\xi_{xx}+\xi_{yy}=0\\
\eta_{xx}+\eta_{yy}=0.
\end{array}\right.$$\\
By solving the two-dimensional Laplace equation, we can get a large
number of solutions to the system (3.3). We first look at one simple
case
\\$$\left\{\begin{array}{ll}
\tau=c_1+c_2t\\
\xi=c_3+c_4x+c_5y\\
\eta=c_6-c_5x+c_4y\\
\phi=c_2-2c_4
\end{array}\right.\eqno(3.4)$$\\
Thus we get
\\$$\left\{\begin{array}{ll}
{\bf v_1}=\frac{\partial}{\partial t}\\
{\bf v_2}=\frac{\partial}{\partial x}\\
{\bf v_3}=\frac{\partial}{\partial y}\\
{\bf v_4}=t\frac{\partial}{\partial t}+\frac{\partial}{\partial w}\\
{\bf v_5}=y\frac{\partial}{\partial x}-x\frac{\partial}{\partial y}\\
{\bf v_6}=x\frac{\partial}{\partial x}+y\frac{\partial}{\partial y}-2\frac{\partial}{\partial w}\\
\end{array}\right.\eqno(3.5)$$\\
The corresponding one-parameter transformation groups are
\\$$\left\{\begin{array}{ll}
G_1:(x,y,t+\varepsilon,w)\\
G_2:(x+\varepsilon,y,t,w)\\
G_3:(x,y+\varepsilon,t,w)\\
G_4:(x,y,e^{\varepsilon}t,w+\varepsilon)\\
G_5:(x+\varepsilon y,y-\varepsilon x,t,w)\\
G_6:(e^{\varepsilon}x,e^{\varepsilon}y,t,w-2\varepsilon)\\
\end{array}\right.$$\\
Equivalently, if $w=f(x,y,t)$ is a solution to (3.2), then the
following are also solutions to (3.2):
\\$$\left\{\begin{array}{ll}
w^{(1)}=f(x,y,t-\varepsilon)\\
w^{(2)}=f(x-\varepsilon,y,t)\\
w^{(3)}=f(x,y-\varepsilon,t)\\
w^{(4)}=f(x,y,e^{-\varepsilon}t)+\varepsilon\\
w^{(5)}=f(x-\varepsilon y,y+\varepsilon x,(1+\varepsilon^2)t)\\
w^{(6)}=f(e^{-\varepsilon}x,e^{-\varepsilon}y,t)-2\varepsilon\\
\end{array}\right.$$\\

For example we examine $w^{(6)}$,
\\$$\begin{array}{ll}
&(e^{w^{(6)}})_t-w^{(6)}_{xx}-w^{(6)}_{yy}\cr
&=(e^{-2\varepsilon}e^f)_t-e^{-2\varepsilon}f_{xx}-e^{-2\varepsilon}f_{yy}\cr
&=e^{-2\varepsilon}((e^f)_t-f_{xx}-f_{yy})\cr &=0.
\end{array}$$\\

Next we have the following structure constants table such that the
entry in $i$-row and $j$-volume represents $[{\bf v_i},{\bf v_j}]$:
\\$$\begin{tabular}{c|cccccc}
\hline
Lie&${\bf v_1}$&${\bf v_2}$&${\bf v_3}$&${\bf v_4}$&${\bf
v_5}$&${\bf v_6}$\\
\hline
${\bf v_1}$&$0$&$0$&$0$&${\bf v_1}$&$0$&$0$\\
${\bf v_2}$&$0$&$0$&$0$&$0$&$-{\bf v_3}$&${\bf v_2}$\\
${\bf v_3}$&$0$&$0$&$0$&$0$&${\bf v_2}$&${\bf v_3}$\\
${\bf v_4}$&$-{\bf v_1}$&$0$&$0$&$0$&$0$&$0$\\
${\bf v_5}$&$0$&${\bf v_3}$&$-{\bf v_2}$&$0$&$0$&$0$\\
${\bf v_6}$&$0$&$-{\bf v_2}$&$-{\bf v_3}$&$0$&$0$&$0$\\
\hline
\end{tabular}\eqno(3.6)$$\\

Using the formula
\\$$\begin{array}{ll}
Ad(\exp(\varepsilon{\bf v_i})){\bf
v_j}&=\sum_{n=0}^{\infty}\frac{\varepsilon^n}{n!}(ad {\bf
v_i})^n({\bf v_j})\cr &={\bf v_j}-\varepsilon[{\bf v_i},{\bf
v_j}]+\frac{\varepsilon^2}{2}[{\bf v_i},[{\bf v_i},{\bf
v_j}]]-\cdot\cdot\cdot.
\end{array}$$\\
We get the adjoint representation table for (3.5):
\\$$\begin{tabular}{c|cccccc}
\hline Ad&${\bf v_1}$&${\bf v_2}$&${\bf v_3}$&${\bf v_4}$&${\bf
v_5}$&${\bf v_6}$\\
\hline
${\bf v_1}$&${\bf v_1}$&${\bf v_2}$&${\bf v_3}$&${\bf v_4}-\varepsilon{\bf v_1}$&${\bf v_5}$&${\bf v_6}$\\
${\bf v_2}$&${\bf v_1}$&${\bf v_2}$&${\bf v_3}$&${\bf v_4}$&${\bf v_5}+\varepsilon{\bf v_3}$&${\bf v_6}-\varepsilon{\bf v_2}$\\
${\bf v_3}$&${\bf v_1}$&${\bf v_2}$&${\bf v_3}$&${\bf v_4}$&${\bf v_5}-\varepsilon{\bf v_2}$&${\bf v_6}-\varepsilon{\bf v_3}$\\
${\bf v_4}$&$e^{\varepsilon}{\bf v_1}$&${\bf v_2}$&${\bf v_3}$&${\bf v_4}$&${\bf v_5}$&${\bf v_6}$\\
${\bf v_5}$&${\bf v_1}$&$\cos(\varepsilon){\bf v_2}-\sin(\varepsilon){\bf v_3}$&
$\cos(\varepsilon){\bf v_3}+\sin(\varepsilon){\bf v_2}$&${\bf v_4}$&${\bf v_5}$&${\bf v_6}$\\
${\bf v_6}$&${\bf v_1}$&$e^{\varepsilon}{\bf v_2}$&$e^{\varepsilon}{\bf v_3}$&${\bf v_4}$&${\bf v_5}$&${\bf v_6}$\\
\hline
\end{tabular}\eqno(3.7)$$\\

Now we use the adjoint representation table to give the
classification of subalgebras of (3.5). Given a vector
\\$${\bf v}=a_1{\bf v_1}+a_2{\bf v_2}+a_3{\bf v_3}+a_4{\bf v_4}+a_5{\bf v_5}+a_6{\bf v_6},$$\\
we first assume $a_6\neq0$, so after scaling, we can make $a_6=1$:
\\$${\bf v}=a_1{\bf v_1}+a_2{\bf v_2}+a_3{\bf v_3}+a_4{\bf v_4}+a_5{\bf v_5}+{\bf v_6}.$$\\
If we act on ${\bf v}$ by $Ad(\exp((a_2-a_5(a_3+a_2a_5)){\bf v_2}))
$ and $Ad(\exp((a_3+a_2a_5){\bf v_3})) $ respectively, we can make
the coefficients of ${\bf v_2}$ and ${\bf v_3}$ vanish:
\\$${\bf v^{(1)}}=Ad(\exp((a_3+a_2a_5){\bf v_3}))\circ Ad(\exp((a_2-a_5(a_3+a_2a_5)){\bf v_2})){\bf v}
=a_1{\bf v_1}+a_4{\bf v_4}+a_5{\bf v_5}+{\bf v_6}.$$\\
Next we act on ${\bf v^{(1)}}$ by $Ad(\exp(a_1{\bf v_1})$ to cancel
the the coefficient of ${\bf v_1}$, so finally ${\bf v}$ is
equivalent to ${\bf v^{(2)}}=a_4{\bf v_4}+a_5{\bf v_5}+a_6$ under
the adjoint representation. In other words, every one-dimensional
subalgebra generated by ${\bf v}$ with $a_6\neq0$ is equivalent to
the subalgebra spanned by $a_4{\bf v_4}+a_5{\bf v_5}+a_6$.

The remaining one-dimensional subalgebras are spanned by vector with
$a_6=0$. If $a_5\neq0$, by scaling we make $a_5=1$, and then act on
${\bf v}$ by $Ad(\exp(-a_3{\bf v_2}))$ and $Ad(\exp(a_2{\bf v_3}))$
respectively so that ${\bf v}$ is equivalent to ${\bf
v^{(1)}}=a_4{\bf v_4}+{\bf v_5}$.

Next, if $a_5=a_6=0$ and $a_4\neq0$, consider ${\bf v}=a_1{\bf
v_1}+a_2{\bf v_2}+a_3{\bf v_3}+{\bf v_4}$. First act on it by
$Ad(\exp(a_1{\bf v_1}))$,
\\$${\bf v^{(1)}}=Ad(\exp(a_1{\bf v_1})){\bf v}=a_2{\bf v_2}+a_3{\bf v_3}+{\bf v_4}.$$\\
If $a_2=0$, then ${\bf v^{(1)}}=a_3{\bf v_3}+{\bf v_4}$. Otherwise,
\\$${\bf v^{(2)}}=Ad(\exp(\arctan(a_3/a_2){\bf v_5})){\bf v^{(1)}}=l{\bf v_2}+{\bf v_4},$$\\
where $l=l(a_2,a_3)$. Further we can use $ $ and $ $ act on ${\bf
v^{1}}$ and ${\bf v^{(2)}}$ respectively to scale $a_3$ and $l$.
Thus any one-dimensional subalgebra spanned by ${\bf v}$ with
$a_5=a_6=0$ and $a_4\neq0$ is equivalent to the subalgebra spanned
by either ${\bf v_4}$, ${\bf v_4}+{\bf v_2}$, ${\bf v_4}-{\bf v_2}$,
${\bf v_4}+{\bf v_3}$ or ${\bf v_4}-{\bf v_3}$.

If $a_4=a_5=a_6=0$ and $a_1\neq0$, let ${\bf v}={\bf v_1}+a_2{\bf
v_2}+a_3{\bf v_3}$. If $a_2=0$, then ${\bf v}={\bf v_1}+a_3{\bf
v_3}$. Otherwise
\\$${\bf v^{(1)}}=Ad(\exp(\arctan(a_3/a_2){\bf v_5})){\bf v}={bf v_1}+l{\bf v_2},$$\\
where $l=l(a_2,a_3)$. By further scaling $a_3$ and $l$ using the
adjoint representation, we finally get the result that when
$a_4=a_5=a_6=0$ and $a_1\neq0$, the subalgebra spanned by {\bf v} is
equivalent to either ${\bf v_1}$, ${\bf v_1}+{\bf v_2}$, ${\bf
v_1}-{\bf v_2}$, ${\bf v_1}+{\bf v_3}$, ${\bf v_1}-{\bf v_3}$.

If $a_1=a_4=a_5=a_6=0$, then the subalgebra spanned by ${\bf v}$ is
equivalent to either ${\bf v_2}+a_3{\bf v_3}$ or ${\bf v_3}$.

If we further allow the discrete symmetry for example
$(\xi,\eta,\tau,\phi)\mapsto(-\xi,\eta,\tau,\phi)$, then ${\bf
v_1}-{\bf v_2}$ is mapped to ${\bf v_1}+{\bf v_2}$. So the following
theorem holds:\\

{\bf Theorem 3.1.} {\it the operators in (3.5) generate an optimal
system $\mathcal {S}$
\begin{enumerate}
\item[${\bf (a)}$]${\bf v_6}+a_4{\bf v_4}+a_5{\bf
v_5}$, $a_6\neq0$;
\item[${\bf (b)}$]${\bf v_5}+a_4{\bf v_4}$, $a_6=0$,
$a_5\neq0$;
\item[${\bf (c_1)}$]${\bf v_4}$, $a_5=a_6=0$, $a_4\neq0$;
\item[${\bf (c_2)}$]${\bf v_4}+{\bf v_2}$, $a_5=a_6=0$, $a_4\neq0$;
\item[${\bf (c_3)}$]${\bf v_4}+{\bf v_3}$, $a_5=a_6=0$, $a_4\neq0$;
\item[${\bf (d_1)}$]${\bf v_1}$, $a_4=a_5=a_6=0$, $a_1\neq0$;
\item[${\bf (d_2)}$]${\bf v_1}+{\bf v_2}$, $a_4=a_5=a_6=0$, $a_1\neq0$;
\item[${\bf (d_3)}$]${\bf v_1}-{\bf v_3}$, $a_4=a_5=a_6=0$, $a_1\neq0$;
\item[${\bf (e)}$]${\bf v_2}+a_3{\bf v_3}$, $a_1=a_4=a_5=a_6=0$, $a_2\neq0$;
\item[${\bf (f)}$]${\bf v_3}$, $a_1=a_2=a_4=a_5=a_6=0$.\\
\end{enumerate}
}

We calculate two examples to show how to use the subalgebras to find
symmetric solutions of (3.2) or (3.1).

From ${\bf v}={\bf v_4}+{\bf v_2}=t\frac{\partial}{\partial
t}+\frac{\partial}{\partial x}+\frac{\partial}{\partial w}$, its
characteristics are derived from
\\$$\frac{dt}{t}=\frac{dx}{1}=\frac{dw}{1}.$$\\
So let $z=e^y/t$ and $e^w=t\omega(z)$, then we have
\\$$z^2\omega\omega''-z^2(\omega')^2+z\omega^2\omega'-\omega^3=0.\eqno(3.8)$$\\

From ${\bf v}=2{\bf v_4}+{\bf v_5}=2t\frac{\partial}{\partial
t}+x\frac{\partial}{\partial x}+y\frac{\partial}{\partial y}$, let
$\varepsilon=\frac{x}{\sqrt t}$, $\eta=\frac{y}{\sqrt t}$, and
$e^w=\omega(\varepsilon,\eta)$, then we have
\\$$\omega_{\varepsilon\varepsilon}+\omega_{\eta\eta}=\frac{\omega_{\varepsilon}^2
+\omega_{\eta}^2}{\omega}-\frac{1}{2}\omega(\varepsilon\omega_{\varepsilon}+\eta\omega_{\eta}).\eqno(3.9)$$\\

Note that we only used a quite simple solution (3.4) of (3.3) to
analyze the symmetries of (3.2). If we use other solutions of (3.3),
we expect to find many more symmetries of (3.2). The whole question
lies in finding solutions for the two-dimensional Laplace equation.
In contrast to the linear heat equation, amazingly (3.1) has various
symmetries. We believe study of (3.3) will lead to some significant
results for Ricci flow on surfaces. For example, we can consider
\\$$\left\{\begin{array}{ll}
\tau=c_1+c_2t\\
\xi=c_3(x^2-y^2)+c_4xy+c_5\\
\eta=\frac{1}{2}c_4(y^2-x^2)+2c_3xy+c_6\\
\phi=c_2-4c_3x-2c_4y,
\end{array}\right.\eqno(3.10)$$\\
or more complicated case
\\$$\left\{\begin{array}{ll}
\tau=c_1+c_2t\\
\xi=(c_3\cos x+c_4\sin x)e^y+(c_4\cos y+c_6\sin y)e^x+c_7\\
\eta=(c_4\cos x-c_3\sin x)e^y+(-c_6\cos y+c_5\sin y)e^x+c_8\\
\phi=c_2+2c_3e^y\sin x-2c_4e^y\cos x-2c_5e^x\cos y-2c_6e^x\sin y.
\end{array}\right.\eqno(3.11)$$\\

\section{Hyperbolic geometric flow on Riemann surfaces}

In this section, we consider the hyperbolic geometric flow on
Riemann surfaces. By the same statements as in the beginning of
section 3, the hyperbolic geometric flow
\\$$\frac{\partial^2}{\partial t^2}g=-2Rc$$\\
can be simplified to
\\$$u_{tt}=\triangle\ln u,\eqno(4.1)$$\\
Let $w=\ln u$, then
\\$$\triangle(x,y,t,w^{(2)})=e^ww_{tt}+e^ww_t^2-w_{xx}-w_{yy}=0.\eqno(4.2)$$\\
The following initial problem has been studied in [10],
\\$$\left\{\begin{array}{ll}
u_{tt}-(\ln u)_{xx}=0\\
t=0:\quad u=u_0(x),\quad u_t=u_t(x).
\end{array}\right.$$\\
Given any initial metric only dependenting on one space variable, if
the initial velocity is large enough, then the solution of (4.1)
exists for all positive time and the scalar curvature is uniformly
bounded. Otherwise, the solution exists only finite time and the
scalar curvature goes to infinity as $t$ goes to the maximal time.

The Jacobian matrix of (4.2) is
\\$$\mathcal {J}_{\triangle}(x,y,t;w;w_x,w_y,w_t;w_{xx},w_{xy},w_{xt},w_{yy},w_{yt},w_{tt})
=(0,0,0;e^w(w_{tt}+w_t^2);0,0,2e^ww_t;-1,0,0,-1,0,e^w).$$\\
Thus the original equation (4.2) is of maximal rank everywhere in
$\mathscr{S}_{\triangle}$. Obviously (4.2) is locally solvable. So
(4.2) is nondegenerate and we can apply Theorem 2.6.

Given a vector Given a vector
\\$${\bf v}=\xi(x,y,t,w)\frac{\partial}{\partial x}+\eta(x,y,t,w)\frac{\partial}{\partial y}
+\tau(x,y,t,w)\frac{\partial}{\partial t}+\phi(x,y,t,w)\frac{\partial}{\partial w},$$\\
we have
\\$$pr^{(2)}{\bf
v}(\triangle(x,y,t,w^{(2)}))=\phi e^w(w_{tt}+w_t^2)+2\phi^te^ww_t+\phi^{tt}e^w-\phi^{xx}-\phi^{yy}.$$\\
We apply the formulas for $\phi^t$, $\phi^{tt}$, $\phi^{xx}$,
$\phi^{yy}$ derived in section 2, since $pr^{(2)}{\bf
v}(\triangle(x,y,t,w^{(2)}))=0$ whenever (4.2) holds, we use
$w_{yy}=e^ww_{tt}+e^ww_t^2-w_{xx}$ to cancel $w_{yy}$ and set
coefficients of every monomial zero. The coefficient table is
\\$$\begin{tabular}{rr|rr}
\hline
monomial&coefficient&monomial&coefficient\\
$e^ww_{tt}:$&$\phi-2\tau_t+2\eta_y=0$&$e^ww_t^2:$&$\phi+\phi_w+\phi_{ww}-2\tau_{tw}-2\tau_t+2\eta_y=0$\\
$e^ww_t:$&$2\phi_t+2\phi_{tw}-\tau_{tt}=0$&$e^ww_tw_x:$&$-2\xi_t-2\xi_{tw}=0$\\
$e^ww_tw_y:$&$-2\eta_t-2\eta_{tw}=0$&$e^ww_xw_t^2:$&$-\xi_w-\xi_{ww}=0$\\
$e^ww_yw_t^2:$&$\eta_w-\eta_{ww}=0$&$e^ww_t^3:$&$-\tau_w-\tau_{ww}=0$\\
$e^w:$&$\phi_{tt}=0$&$e^ww_y:$&$-\eta_{tt}=0$\\
$e^ww_x:$&$-\xi_{tt}=0$&$e^ww_{xt}:$&$-2\xi_t=0$\\
$e^ww_{yt}:$&$-2\eta_t=0$&$e^ww_tw_{tt}:$&$-2\tau_w=0$\\
$e^ww_yw_{tt}:$&$2\eta_w=0$&$e^ww_tw_{yt}:$&$-2\eta_w=0$\\
$e^ww_tw_{xt}:$&$-2\xi_w=0$&$1:$&$-\phi_{xx}-\phi_{yy}=0$\\
$w_y:$&$-2\phi_{yw}+\eta_{xx}+\eta_{yy}=0$&$w_x:$&$-2\phi_{xw}+\xi_{xx}+\xi_{yy}=0$\\
$w_t:$&$\tau_{xx}+\tau_{yy}=0$&$w_y^2:$&$-\phi_{ww}+2\eta_{yw}=0$\\
$w_xw_y:$&$2\xi_{yw}+2\eta_{xw}=0$&$w_yw_t:$&$2\tau_{yw}=0$\\
$w_y^3:$&$\eta_{ww}=0$&$w_y^2w_x:$&$\xi_{ww}=0$\\
$w_y^2w_t:$&$\tau_{ww}=0$&$w_{xx}:$&$2\xi_x-2\eta_y=0$\\
$w_{yt}:$&$2\tau_y=0$&$w_{xy}:$&$2\xi_y+2\eta_x=0$\\
$w_yw_{xx}:$&$-2\eta_w=0$&$w_xw_{xx}:$&$2\xi_w=0$\\
$w_yw_{xy}:$&$-2\xi_w=0$&$w_yw_{yt}:$&$-2\tau_w=0$\\
$w_x^2:$&$-\phi_{ww}+2\xi_{xw}=0$&$w_xw_t:$&$-2\tau_{xw}=0$\\
$w_x^3:$&$\xi_{ww}=0$&$w_x^2w_y:$&$\eta_{ww}=0$\\
$w_x^2w_t:$&$\tau_{ww}=0$&$w_{xt}:$&$2\tau_x=0$\\
$w_xw_{xy}:$&$2\eta_w=0$&$w_xw_{xt}:$&$2\tau_w=0$\\
\hline
\end{tabular}$$\\

Finally, we get
\\$$\left\{\begin{array}{ll}
\xi=\xi(x,y)\\
\eta=\eta(x,y)\\
\tau=c_1+c_2t\\
\phi=2c_2-2\xi_x\\
\xi_x-\eta_y=0\\
\eta_x+\xi_y=0
\end{array}\right.\eqno(4.3)$$\\
This is the same as (3.3) except the coefficient $2$ of $c_2$ in
$\phi$. This is due to the fact that in wave type equation (4.2), we
differentiate $w$ twice along the time variable $t$.

If we choose
\\$$\left\{\begin{array}{ll}
\tau=c_1+c_2t\\
\xi=c_3+c_4x+c_5y\\
\eta=c_6-c_5x+c_4y\\
\phi=2c_2-2c_4
\end{array}\right.\eqno(4.4)$$\\
Then we get
\\$$\left\{\begin{array}{ll}
{\bf v_1}=\frac{\partial}{\partial t}\\
{\bf v_2}=\frac{\partial}{\partial x}\\
{\bf v_3}=\frac{\partial}{\partial y}\\
{\bf v_4}=t\frac{\partial}{\partial t}+2\frac{\partial}{\partial w}\\
{\bf v_5}=y\frac{\partial}{\partial x}-x\frac{\partial}{\partial y}\\
{\bf v_6}=x\frac{\partial}{\partial x}+y\frac{\partial}{\partial y}-2\frac{\partial}{\partial w}\\
\end{array}\right.\eqno(4.5)$$\\
The corresponding one-parameter transformation groups are
\\$$\left\{\begin{array}{ll}
G_1:(x,y,t+\varepsilon,w)\\
G_2:(x+\varepsilon,y,t,w)\\
G_3:(x,y+\varepsilon,t,w)\\
G_4:(x,y,e^{\varepsilon}t,w+2\varepsilon)\\
G_5:(x+\varepsilon y,y-\varepsilon x,t,w)\\
G_6:(e^{\varepsilon}x,e^{\varepsilon}y,t,w-2\varepsilon)\\
\end{array}\right.$$\\
Equivalently, if $w=f(x,y,t)$ is a solution to (4.2), then the
following are also solutions to (4.2):
\\$$\left\{\begin{array}{ll}
w^{(1)}=f(x,y,t-\varepsilon)\\
w^{(2)}=f(x-\varepsilon,y,t)\\
w^{(3)}=f(x,y-\varepsilon,t)\\
w^{(4)}=f(x,y,e^{-\varepsilon}t)+2\varepsilon\\
w^{(5)}=f(x-\varepsilon y,y+\varepsilon x,\sqrt{1+\varepsilon^2}t)\\
w^{(6)}=f(e^{-\varepsilon}x,e^{-\varepsilon}y,t)-2\varepsilon\\
\end{array}\right.$$\\

The the commutator table and the adjoint representation of (4.5) are
the same as those of (3.5), namely (3.6) and (3.7), so we omit them
here.

Similarly, we have\\

{\bf Theorem 4.1.} {\it the operators in (4.5) generate an optimal
system $\mathcal {S}$
\begin{enumerate}
\item[${\bf (a)}$]${\bf v_6}+a_4{\bf v_4}+a_5{\bf
v_5}$, $a_6\neq0$;
\item[${\bf (b)}$]${\bf v_5}+a_4{\bf v_4}$, $a_6=0$,
$a_5\neq0$;
\item[${\bf (c_1)}$]${\bf v_4}$, $a_5=a_6=0$, $a_4\neq0$;
\item[${\bf (c_2)}$]${\bf v_4}+{\bf v_2}$, $a_5=a_6=0$, $a_4\neq0$;
\item[${\bf (c_3)}$]${\bf v_4}+{\bf v_3}$, $a_5=a_6=0$, $a_4\neq0$;
\item[${\bf (d_1)}$]${\bf v_1}$, $a_4=a_5=a_6=0$, $a_1\neq0$;
\item[${\bf (d_2)}$]${\bf v_1}+{\bf v_2}$, $a_4=a_5=a_6=0$, $a_1\neq0$;
\item[${\bf (d_3)}$]${\bf v_1}-{\bf v_3}$, $a_4=a_5=a_6=0$, $a_1\neq0$;
\item[${\bf (e)}$]${\bf v_2}+a_3{\bf v_3}$, $a_1=a_4=a_5=a_6=0$, $a_2\neq0$;
\item[${\bf (f)}$]${\bf v_3}$, $a_1=a_2=a_4=a_5=a_6=0$.\\
\end{enumerate}
}

From ${\bf v}={\bf v_4}+{\bf v_2}=t\frac{\partial}{\partial
t}+\frac{\partial}{\partial x}+2\frac{\partial}{\partial w}$, its
characteristics are derived from
\\$$\frac{dt}{t}=\frac{dx}{1}=\frac{dw}{2}.$$\\
So let $z=e^y/t$ and $e^w=t^2\omega(z)$, then we have
\\$$z^2(\omega^2-\omega)\omega''-z(2\omega^2+\omega)\omega'+z^2(\omega')^2+2\omega^3=0.\eqno(4.6)$$\\

From ${\bf v}={\bf v_4}+{\bf v_5}=t\frac{\partial}{\partial
t}+x\frac{\partial}{\partial x}+y\frac{\partial}{\partial y}$, let
$\varepsilon=\frac{x}{t}$, $\eta=\frac{y}{t}$, and
$e^w=\omega(\varepsilon,\eta)$, then we have
\\$$(\varepsilon^2\omega^2-\omega)\omega_{\varepsilon\varepsilon}
+(\eta^2\omega^2-\omega)\omega_{\eta\eta}+2\varepsilon\eta\omega^2\omega_{\varepsilon\eta}
+\omega_{\varepsilon}^2+\omega_{\eta}^2=0.\eqno(4.7)$$\\

\section{Warped products on ${\bf\mathcal {S}^{n+1}}$ of both flows}

In contrast to the Ricci flow and hyperbolic geometric flow on
surfaces, the warped products on $\mathcal {S}^{n+1}$ of both flows
admit few symmetries. We omit the somewhat laboring calculations and
only state the results here.

Recall that warped product on $\mathcal {S}^{n+1}$ is of the form
\\$$g=\varphi^2(x,t)dx^2+\psi^2(x,t)g_{can},$$\\
where $g_{can}$ denotes the canonical metric on $\mathcal {S}^n$. By
a change of coordinate
\\$$s(x)=\int^x_0\varphi(x)dx,\eqno(5.1)$$\\
the evolutions of $\varphi(s,t)$ and $\psi(s,t)$ under Ricci flow
and hyperbolic geometric flow are the followings respectively:
\\$$\left\{\begin{array}{ll}
\varphi_t=n\frac{\psi_{ss}}{\psi}\varphi\\
\psi_t=\psi_{ss}-(n-1)\frac{1-\psi^2_s}{\psi}
\end{array}\right.\eqno(5.2)$$\\
under Ricci flow, and
\\$$\left\{\begin{array}{ll}
\varphi_{tt}=n\frac{\psi_{ss}}{\psi}\varphi-\frac{\varphi^2_t}{\varphi}\\
\psi_{tt}=\psi_{ss}-(n-1)\frac{1-\psi^2_s}{\psi}-\frac{\psi^2_t}{\psi}
\end{array}\right.\eqno(5.3)$$\\
under hyperbolic geometric flow. Due to the change of coordinate
(5.1), we can only consider the second equations of (5.2) and (5.3),
namely
\\$$\psi_t=\psi_{ss}-(n-1)\frac{1-\psi^2_s}{\psi},\eqno(5.4)$$\\
and
\\$$\psi_{tt}=\psi_{ss}-(n-1)\frac{1-\psi^2_s}{\psi}-\frac{\psi^2_t}{\psi}.\eqno(5.5)$$\\

Given a vector
\\$${\bf v}=\xi(s,t,\psi)\frac{\partial}{\partial s}+\tau(s,t,\psi)\frac{\partial}{\partial t}
+\phi(s,t,\psi)\frac{\partial}{\partial\psi}$$\\
We have the following two theorems:\\

{\bf Theorem 5.1(symmetries for warped product of Ricci flow).} {\it
For equation (5.4), we have\\
${\bf when\quad n=1:}$ $$ \left\{\begin{array}{ll}
{\bf v_1}=\frac{\partial}{\partial s}\\
{\bf v_2}=\frac{\partial}{\partial t}\\
{\bf v_3}=s\frac{\partial}{\partial s}+t\frac{\partial}{\partial t}
\end{array}\right.$$\\
These are translations and dilatation.\\
${\bf when\quad n=2:}$ $$ \left\{\begin{array}{ll}
{\bf v_1}=\frac{\partial}{\partial s}\\
{\bf v_2}=\frac{\partial}{\partial t}\\
{\bf v_3}=t\frac{\partial}{\partial s}+s\frac{\partial}{\partial t}
\end{array}\right.$$\\
These are translations and hyperbolic rotation.\\
${\bf when\quad n>2:}$ $$ \left\{\begin{array}{ll}
{\bf v_1}=\frac{\partial}{\partial s}\\
{\bf v_2}=\frac{\partial}{\partial t}\\
\end{array}\right.$$
These are only translations.\\}

{\bf Theorem 5.2(symmetries for warped product of hyperbolic
geometric flow).} {\it For equation (5.5), we have\\
${\bf when\quad n=1:}$ the equation becomes linear heat equation, so
\\$$
\left\{\begin{array}{ll}
{\bf v_1}=\frac{\partial}{\partial s}\\
{\bf v_2}=\frac{\partial}{\partial t}\\
{\bf v_3}=\psi\frac{\partial}{\partial\psi}\\
{\bf v_4}=s\frac{\partial}{\partial s}+2t\frac{\partial}{\partial t}\\
{\bf v_5}=2t\frac{\partial}{\partial
s}-s\psi\frac{\partial}{\partial\psi}\\
{\bf v_6}=4ts\frac{\partial}{\partial
s}+4t^2\frac{\partial}{\partial
t}-(s^2+2t)\psi\frac{\partial}{\partial\psi},
\end{array}\right.$$\\
and the infinite-dimensional subalgebra
\\$${\bf v_{\alpha}}=\alpha(s,t)\frac{\partial}{\partial\psi},$$\\
where $\alpha$ is an arbitrary solution of the heat equation.\\
${\bf when\quad n>1:}$ $$ \left\{\begin{array}{ll}
{\bf v_1}=\frac{\partial}{\partial s}\\
{\bf v_2}=\frac{\partial}{\partial t}\\
\end{array}\right.$$\\}

\section{Further discussion}

Ricci flow is a powerful tool to understand the geometry and
topology of Riemann manifolds. Any symmetry and exact solution of
its equation will help us understand its behavior for general cases
and the singularity formation, further the basic topological and
geometrical properties as well as analytic properties of the
underlying manifolds. The hyperbolic geometric flow is the
hyperbolic version of Ricci flow. It is also closely related to the
Einstein equation. Any symmetry and exact solution of it can help us
find new solutions of the Einstein Equation which plays significant
role in general relativity and modern theoretical physics.

The techniques we use in this paper, namely the theory of
group-invariant solutions for differential equations is a powerful
tool to analyze various differential equations. In fact, it can also
be used in normal functions and systems of differential equations,
see \cite{O}. We hope in the future this method can be generalized
to tensor equations so that we can use it to analyze complicated
systems on manifolds. This theory of group-invariant solutions is an
application of Lie groups to differential equations. Lie group plays
a fundamental role in modern mathematics since it has significant
influence on almost all the branches of mathematics.\\

\section{Appendix: Derivation of the evolutions of warped products on
${\bf\mathcal {S}^{n+1}}$}

In this appendix, we derive evolution equations from the metrics on
$\mathcal {S}^{n+1}$:
\\$$g = \varphi(x)^2dx\otimes dx + \psi(x)^2ds^2_n$$\\
where $ds^2_n$ denotes the canonical metric of constant curvature 1
on $\mathcal {S}^n$.

By introducing a new coordinate the distance $s$ to the equator
given by
\\$$s(x) = \int_0^x\varphi(x)dx,$$\\
the metrics can be simplified to
\\$$g = ds^2 + \psi(s)^2ds^2_n = ds^2 + g_s.\eqno(7.1)$$\\
In the following we will consider (7.1).

Metrics (7.1) are standard warped products which are studied in
details in \cite{PP}. Our discussion follows closely with Chapter 3
of
that book.\\

Given a Riemannian manifold $(\mathcal{M},g)$, first we define the
{\it Hessian} function of $f$ as a symmetric $(0,2)$-tensor
\\$$Hess~f(X,Y) = \nabla_X\nabla_Yf-\nabla_{\nabla_XY}f,$$\\
where $X, Y$ are vector fields on $\mathcal{M}$. Thus
\\$$Hess~f(X,Y) = g(\nabla_X(\nabla f),Y).$$\\

Second we will say that $s: U\rightarrow\mathbb{R}$, where $U\subset
(\mathcal{M},g)$ is open, is a {\it distance function} if $|\nabla
s|\equiv1$ on $U$. We shall use the following {\it Gauss equation}
about distance function
\\$$g(R(X,Y)Z,W)=g_s(R^s(X,Y)Z,W)-II(Y,Z)II(X,W)+II(X,Z)II(Y,W),$$\\
here $X,Y,Z,W$ are tangent to the level sets $U_s$ and
\\$$II(U,V)=Hess~s(U,V)$$\\
is the classical second fundamental form.

Particularly in our case, for the rotationally symmetric metrics
(7.1),
\\$$\begin{array}{ll}
  2Hess~s &= L_{\partial s}g_s\cr
          &= L_{\partial s}(\psi^2ds^2_n)\cr
          &= \partial_s(\psi^2)ds^2_n+\psi^2L_{\partial
             s}(ds^2_n)\cr
          &= 2\psi(\partial_s\psi)ds^2_n\cr
          &= 2\frac{\partial_s\psi}{\psi}g_s.
    \end{array}$$\\
Thus
\\$$Hess~s = \frac{\partial_s\psi}{\psi}g_s.$$\\

Using that $g_s$ is the metric of curvature $\frac{1}{\psi^2}$ on
the sphere, we get
\\$$g_s(R^s(X,Y)V,W) = \frac{1}{\psi^2}g_s(X\wedge Y,W\wedge V).$$\\

Combining this with $II=Hess~s$, from the Guass equation we obtain
\\$$
  g(R(X,Y)V,W) = \frac{1-(\partial_s\psi)^2}{\psi^2}g_s(X\wedge Y,W\wedge
  V).\eqno(7.2)$$\\

From another important formula
\\$$(\nabla_{\partial_s}Hess~s)(X,Y)+Hess^2~s(X,Y)=-R(X,\partial_s,\partial_s,Y),$$\\
and in our case
\\$$\begin{array}{ll}
     \nabla_{\partial_s}Hess~s
  &= \nabla_{\partial_s}(\frac{\partial_s\psi}{\psi}g_s)\cr
  &= \partial_s(\frac{\partial_s\psi}{\psi})g_s+\frac{\partial_s\psi}{\psi}\nabla_{\partial_s}(g_s)\cr
  &= \frac{(\partial^2_s\psi)\psi-(\partial_s\psi)^2}{\psi^2}g_s\cr
  &= \frac{\partial^2_s\psi}{\psi}g_s-(\frac{\partial_s\psi}{\psi})^2g_s\cr
  &= \frac{\partial^2_s\psi}{\psi}g_s-Hess^2~s,
    \end{array}$$\\
we obtain
\\$$R(\cdot,\partial_s,\partial_s,\cdot) =
-\frac{\partial^2_s\psi}{\psi}.\eqno(7.3)$$\\

(7.2) and (7.3) are just
\\$$K_0=-\frac{\psi_{ss}}{\psi},\qquad K_1=\frac{1-\psi^2_s}{\psi^2},$$\\
where $K_0$ are the sectional curvatures of the $2$-planes
perpendicular to the spheres $\{x\}\times \mathcal {S}^n$, and $K_1$
those of the $2$-planes tangential to these spheres.

Hence for tangential vector $X$,
\\$$\begin{array}{ll}
 Ric(X) &= \sum_{i=1}^{n+1}R(X,E_i)E_i\cr
        &= \sum_{i=1}^nR(X,E_i)E_i+R(X,\partial_s)\partial_s\cr
        &= ((n-1)K_1+K_0)X.
    \end{array}$$\\
And
\\$$Ric(\partial_s)=nK_0.$$\\
Since the metrics (7.1) are Einstein metrics, we write
\\$$Rc=nK_0ds^2+((n-1)K_1+K_0)g_s.$$\\
Finally we obtain
\\$$Rc[g]=-n\frac{\psi_{ss}}{\psi}ds^2+[-\psi\psi_{ss}-(n-1)\psi^2_s+n-1]ds^2_n.\eqno(7.4)$$\\

Differentiate $g$ along $t$, we get
\\$$\frac{\partial}{\partial t}g=\frac{2\varphi_t}{\varphi}ds^2+2\psi\psi_tds^2_n,$$\\
and
\\$$\frac{\partial^2}{\partial t^2}g=\frac{2(\varphi\varphi_{tt}+\varphi^2_t)}{\varphi^2}ds^2
+2(\psi\psi_{tt}+\psi^2_t)ds^2_n.$$\\
So for Ricci flow $\frac{\partial}{\partial t}g=-2Rc$,
\\$$\left\{\begin{array}{ll}
    \frac{\partial}{\partial
    t}\varphi=n\frac{\psi_{ss}}{\psi}\varphi\\
    \frac{\partial}{\partial
    t}\psi=\psi_{ss}-(n-1)\frac{1-\psi^2_s}{\psi}.
    \end{array}\right.$$\\
For hyperbolic geometric flow $\frac{\partial^2}{\partial
t^2}g=-2Rc$,
\\$$\left\{\begin{array}{ll}
    \frac{\partial^2}{\partial
    t^2}\varphi=n\frac{\psi_{ss}}{\psi}\varphi-\frac{\varphi^2_t}{\varphi}\\
    \frac{\partial^2}{\partial
    t^2}\psi=\psi_{ss}-(n-1)\frac{1-\psi^2_s}{\psi}-\frac{\psi^2_t}{\psi}.
    \end{array}\right.$$\\


\begin{thebibliography}{99}

\bibitem{AK} Sigurd B. Angenent and Dan Knopf. {\it An example of
nechpinching for Ricci flow on $\mathcal {S}^{n+1}$}. arXiv.

\bibitem{CZ} Huai-Dong Cao and Xi-Ping Zhu. {\it A complete proof of the
Poincar$\acute{e}$ and geometrization conjectures- application of
the Hamilton-Perelman theory of the Ricci flow}. Asian J. Math.
10(2006), 165-492.

\bibitem{CK} Bennett Chow and Dan Knopf. {\it The Ricci flow: An
introduction}. Mathematical Surveys and Monographs, AMS, Providence,
RI, 2004.

\bibitem{CLN} Bennett Chow, Peng Lu and Lei Ni. {\it Hamilton's Ricci flow}.
Lectures in Contemporary Mathematics, 3, Science Press and Graduate
Studies in Mathematics, 77, AMS, 2006.

\bibitem{CCGG} Bennett Chow, Sun-Chin Chu, David Glickenstein, Christine
Guenther, Jim Isenberg, Tom Ivey, Dan Knopf, Peng Lu, Feng Luo, and
Lei Ni. {\it The Ricci flow: Techniques and Application, Part I:
geometric Aspects}. Mathematical Surveys and Monographs, 135, AMS,
Providence, RI, 2007.

\bibitem{DKL} Wen-Rong Dai, De-Xing Kong and Kefeng Liu. {\it Hyperbolic
geometric flow (I): short-time existence
 and nonlinear stability}. Pure and Applied Mathematics
Quarterly.

\bibitem{H1} Richard S. Hamilton. {\it The Ricci flow on surfaces}.
Mathemat7ics and general relativity (Santa Cruz, CA, 1986), 237-262,
Contemp. Math., 71, AMS, Providence, RI, 1988.

\bibitem{H2} Richard S. Hamilton. {\it The formation of singularities in the
Ricci flow }. Surveys in differential geometry, Vol.II(Cambridge,
MA, 1993), 7-136, Internat. Press, Cambridge, MA, 1995.

\bibitem{KL} De-Xing Kong and Kefeng Liu. {\it Wave character of metrics and
hyperbolic geometric flow}. J. Math. Phys. 48(2007),
103508-1-103508-14.

\bibitem{KLX} De-Xing Kong, Kefeng Liu and De-Liang Xu. {\it The hyperbolic
geometric flow on Riemann surfaces}. Comm. in Partial Differential
Equation 34(2009), 553-580.

\bibitem{O} Peter J. Ovler. {\it Applications of Lie groups to differential
Equations}. GTM, 107. Springer-Verlag, 1993.

\bibitem{PP} Peter Petersen. {\it Riemannian Geometry}. GTM, Springer-Verlag.

\end{thebibliography}
\end{document}